\documentclass[10pt]{amsart}
\usepackage{amsxtra,latexsym,amssymb}
\usepackage{epsfig,rotate,amsthm}

\def\qed{{\hfill $\Box$}}

\def\N{{\mathbb N}}
\def\Z{{\mathbb Z}}
\def\C{{\mathbb C}}

\theoremstyle{theorem}
\newtheorem{thm}{Theorem}[section]
\newtheorem{cor}{Corollary}[section]
\newtheorem{prop}{Proposition}[section]

\newtheorem{lem}{Lemma}[section]
\theoremstyle{definition}
\newtheorem{defn}{Definition}[section]
\theoremstyle{remark}

\newtheorem{rem}{\bf Remark}[section]

\begin{document}

\begin{center}
{\Large {\bf CONSTRUCTING IRREDUCIBLE REPRESENTATIONS OF QUANTUM
    GROUPS $U_{q}(f_{m}(K))$}}\\
\end{center}
\vskip .2in
\begin{center}
Xin Tang 
\bigskip

{\footnotesize Department of Mathematics $\&$ Computer Science, Fayetteville State University,\\
 Fayetteville, NC 28301, U.S.A \hspace{5mm} E-mail:
xtang@uncfsu.edu}

\end{center}

\bigskip

\begin{center}

\begin{minipage}{12cm}
{\bf Abstract}: 
In this paper, we construct families of irreducible representations 
for a class of quantum groups $U_{q}(f_{m}(K))$. First, we give a 
natural construction of irreducible weight representations for 
$U_{q}(f_{m}(K))$ using methods in spectral theory developed by 
Rosenberg. Second, we study the Whittaker model for the center 
of $U_{q}(f_{m}(K))$. As a result, the structure of Whittaker 
representations is determined, and all irreducible Whittaker 
representations are explicitly constructed. Finally, we prove 
that the annihilator of a Whittaker representation is centrally 
generated. 
\medskip


\medskip


\end{minipage}

\end{center}
                         
\section*{ 0. Introduction}
Note that the quantized enveloping algebra $U_{q}(sl_{2})$ has 
played a fundamental role in the study of the quantized enveloping 
algebra $U_{q}(\frak g)$ of any semisimple Lie algebra $\frak g$. In 
order to study the deformations of $U_{q}(\frak g)$, it is natural 
to first study deformations of $U_{q}(sl_{2})$. 

In \cite{JWZ}, a class of algebras $U_{q}(f(K))$ parameterized by 
Laurent polynomials $f(K) \in \C[K,K^{-1}]$ was introduced as 
generalizations of $U_{q}(sl_{2})$. The condition for the existence 
of a Hopf algebra structure on $U_{q}(f(K))$ was determined and 
finite dimensional irreducible representations were explicitly 
constructed. 

Such generalizations yield a family of quantum groups in 
the sense of Drinfeld \cite{D}. For some special parameters 
$f(K)=a(K^{m}-K^{-m}), a \neq 0, m\in \N$, $U_{q}(f(K))$ are 
quantum groups, and all finite dimensional representations 
of $U_{q}(f(K))$ are completely reducible \cite{JWZ}. In 
particular, $U_{q}(f_{m}(K))$ with $f_{m}=\frac{K^{m}-K^{-m}}{q-q^{-1}}$ 
are quantum groups. In this paper, we study the irreducible
representations of these quantum groups $U_{q}(f_{m}(K))$.

As a matter of fact, $U_{q}(f(K))$ can also be realized as Hyperbolic 
algebras \cite{R} (or under the name of Generalized Weyl algebras
\cite{B}). For Hyperbolic algebras, the spectral theory developed 
in \cite{R} is convenient for the explicit construction of irreducible 
weight representations. As an application, we construct all
irreducible weight representations for the quantum groups 
$U_{q}(f_{m}(K))$ via methods developed in \cite{R}.

The Whittaker model for the center $Z(\frak g)$ of the universal 
enveloping algebra $U(\frak g)$ for any semisimple Lie algebra 
$\frak g$ was first studied by Kostant in the seminal paper \cite{K}. 
Whittaker model is closely related to Whittaker equations, and has 
a nice application in the theory of Toda lattice. Via the Whittaker 
model for $Z(\frak g)$, the structure of Whittaker representations 
was determined, and all irreducible Whittaker representations were 
classified in \cite{K}.

The quantum analogue of the Whittaker model for $U_{q}(sl_{2})$ was 
obtained in \cite{O}. For the semisimple Lie algebras of higher ranks, 
the quantum Whittaker model was constructed in \cite{S1} for the
topological version of quantized enveloping algebras via their Coxeter
realizations. In this paper, we try to generalize Kostant's results to 
the quantum groups $U_{q}(f_{m}(K))$. We construct the Whittaker model 
for the center of $U_{q}(f_{m}(K))$. As a result, we determine the 
structure of any Whittaker representation and construct all
irreducible Whittaker representations explicitly. In addition, 
we prove that the annihilator of any Whittaker representation is 
centrally generated.

The above two constructions of irreducible representations 
for $U_{q}(f_{m}(K))$ are perpendicular in the sense that 
Whittaker representations are not weight representations. 
Though the constructions may work for general parameters 
$f(K)$, the calculations are more involved.  

The paper is organized as follows. In Section 1, we recall the
definition and some basic facts about $U_{q}(f(K))$. In Section 2, 
we recall some background about spectral theory and Hyperbolic 
algebras from \cite{R}. Then we illustrate how to realize
$U_{q}(f(K))$ as Hyperbolic algebras, and construct irreducible weight 
representations for $U_{q}(f_{m}(K))$. In Section 3, we describe 
the center of $U_{q}(f_{m}(K))$ and construct the Whittaker model 
for the center. We study the properties of Whittaker representations.

\section{The algebras $U_{q}(f(K))$}
Let $\mathbb{C}$ be the field of complex numbers and $0\neq q \in
\mathbb{C}$ such that $q^{2}\neq 1$. The quantized enveloping algebra 
$U_{q}(sl_{2})$ corresponding to the simple Lie algebra $sl_{2}$ is 
the associative $\C-$algebra generated by $K^{\pm 1},\,E,\, F$ subject to the
following relations:
\[
KE=q^{2}EK,\quad KF=q^{-2}FK,\quad KK^{-1}=K^{-1}K=1;
\]
\[
EF-FE=\frac{K-K^{-1}}{q-q^{-1}}.
\]

It is well-known that $U_{q}(sl_{2})$ is a Hopf algebra with the
following Hopf algebra structure:
\[
\Delta(E)=E\otimes 1+ K\otimes E,\quad \Delta(F)=F\otimes K^{-1}+1\otimes F,\]
\[
\epsilon(E)=0=\epsilon(F),\quad \epsilon(K)=1=\epsilon(K^{-1}),
\]
\[
s(E)=-K^{-1}E,\quad s(F)=-FK,\quad s(K)=K^{-1}.
\]
As generalizations of $U_{q}(sl_{2})$, a class of algebras
$U_{q}(f(K))$ parameterized by Laurent polynomials 
$f(K)\in \mathbb{C}[K,K^{-1}]$ was introduced in \cite{JWZ}. 
We recall their definition here.
\begin{defn}
(See \cite{JWZ}) For any Laurent polynomial $f(K)\in
\mathbb{C}[K,K^{-1}]$, $U_{q}(f(K))$ is defined to be 
the $\mathbb{C}-$algebra generated by $E,\,F,\,K^{\pm 1}$ 
subject to the following relations:
\[
KE=q^{2}EK,\quad KF=q^{-2}FK,
\]
\[
KK^{-1}=K^{-1}K=1,
\]
\[
EF-FE=f(K).
\]
\end{defn}
The ring theoretic properties and finite dimensional representations 
were studied in detail in \cite{JWZ}. We recall some basic results 
without proof.
\begin{prop}
(Prop 3.3 in \cite{JWZ}) Assume $f(K)$ is a non-zero Laurent 
polynomial in $\mathbb{C}[K,\,K^{-1}]$. Then the non-commutative 
algebra $U_{q}(f(K))$ is a Hopf algebra such that 
$K,\,K^{-1}$ are group-like elements, and $E,\, F$ are 
skew primitive elements if and only if $f(K)=a(K^{m}-K^{-m})$ 
with $m=t-s$ and the following conditions are satisfied:
\begin{gather*}
\Delta(K)=K\otimes K,\quad \Delta(K^{-1})=K^{-1}\otimes K^{-1};\\
\Delta(E)=E^{s}\otimes E + E\otimes K^{t},\quad \Delta(F)=K^{-t}\otimes F+F\otimes K^{-s};\\
\epsilon(K)=\epsilon(K^{-1})=1,\quad \epsilon(E)=\epsilon(F)=0;\\
S(K)=K^{-1},\quad S(K^{-1})=K;\\
S(E)=-K^{-s}EK^{-t},\quad S(F)=-K^{t}FK^{s}.
\end{gather*}
\end{prop}
\qed
 
In particular, $U_{q}(f_{m}(K))$ are quantum groups for $f_{m}(K)
=\frac{K^{m}-K^{-m}}{q-q^{-1}}, m\in \mathbb{N}$. When $q$ is not 
a root of unity, the finite dimensional irreducible representations 
were proved to be highest weight and were constructed explicitly in
\cite{JWZ}. Furthermore, we have the following: 
\begin{thm}
(Thm 4.17 in \cite{JWZ}) If $q$ is not a root of unity, then any 
finite dimensional representation $V$ of $U_{q}(f_{m}(K))$ is 
completely reducible.
\end{thm}
\qed

\section{Hyperbolic algebras and their representations}
In this section, we realize $U_{q}(f(K))$ as Hyperbolic algebras and 
use the methods in spectral theory developed in \cite{R} to construct 
irreducible weight representations for $U_{q}(f_{m}(K))$. For the
reader's convenience, we recall some background about spectral theory 
from \cite{R}.

\subsection{Preliminaries on spectral theory}
Spectral theory of abelian categories was started by Gabriel in
\cite{G}. Gabriel defined the injective spectrum of any noetherian 
Grothendieck category with enough injectives. This spectrum consists 
of isomorphism classes of indecomposable injective objects of the
category. Let $R$ be a commutative noetherian ring, then the spectrum 
of the category of all $R-$modules is isomorphic to the prime spectrum 
$Spec(R)$ of $R$ as a scheme. Furthermore, one can reconstruct any 
noetherian commutative scheme $(X,O_{X})$ using the spectrum of the 
category of quasi-coherent sheaves of modules on $X$. The spectrum 
of any abelian category was later on defined by Rosenberg in \cite{R}. 
This spectrum works for any abelian category. Via this spectrum, 
one can reconstruct any quasi-separated and quasi-compact commutative 
scheme $(X,O_{X})$ via the spectrum of the category of quasi-coherent 
sheaves of modules on $X$. 

Though spectral theory is more important for the purpose of 
non-commutative algebraic geometry, it has nice applications 
to representation theory. Note that the spectrum has a natural 
analogue of the Zariski topology and its closed points 
are in a one-to-one correspondence with the irreducible objects 
of the category under certain mild finiteness condition. This 
is the case for the category of representations over an algebra. 
To study irreducible representations, one can study the spectrum 
of the category of all representations, then single out closed 
points of the spectrum with respect to the associated topology. 
 
Let $C_{X}$ be an abelian category and $M,\, N \in C_{X}$ be any two
objects; We say that $M \succ N $ if and only if $N$ is a sub-quotient
of the direct sum of finitely many copies of $M$. It is easy to verify that
$\succ$ is a pre-order. We say $M \approx N$ if and only if $M\succ N$ 
and $N\succ M$. It is obvious that $\approx$ is an equivalence. Let $Spec(X)$ 
be the family of all nonzero objects $M\in C_{X}$ such that for any
non-zero sub-object $N$ of $M$, we have $N\succ M$. 
\begin{defn}
(See \cite{R}) The spectrum of any abelian category is defined to be: 
\[
{\bf Spec(X)}=Spec(X)/\approx.
\]
\end{defn}
 
\subsection{The left spectrum of a ring}
If $C_X$ is the category $A-mod$ of left modules over a ring $A$, then 
it is sometimes convenient to express the points of ${\bf Spec}(X)$ in 
terms of left ideals of the ring $A$. In order to do so, the 
{\it left spectrum} $Spec_l(A)$ was defined in \cite{R}, which 
is by definition the set of all left ideals $p$ of $A$ such 
that $A/p$ is an object of $Spec(X)$. The relation 
$\succ$ on $A-mod$ induces a {\it specialization} relation 
among left ideals, in particular, the specialization relation 
on $Spec_l(A)$. Namely, $A/m\succ A/n$ iff there exists a 
finite subset $x$ of elements of $A$ such that such that the 
ideal $(n:x)=\{a\in A\ |\ ax\subset n\}$ is contained in $m$. 
Following \cite{R}, we denote this by $n\le m$. Note that 
the relation $\le$ is just the inclusion if $n$ is a two-sided 
ideal. In particular, it is the inclusion if the ring $A$ is 
commutative. The map which assigns to an element of $Spec_l(A)$ induces 
a bijection of the quotient $Spec_l(A)/\approx$ of $Spec_l(A)$ 
by the equivalence relation associated with $\le$ onto 
${\bf Spec}(X)$. From now on, we will not distinguish 
$Spec_l(A)/\approx$ from ${\bf Spec}(X)$ and will express 
results in terms of the left spectrum.  

\subsection{Hyperbolic algebra $R\{\xi,\,\theta\}$ and its spectrum}
Hyperbolic algebras are studied by Rosenberg in \cite{R} and by Bavula
under the name of Generalized Weyl algebras in \cite{B}. Hyperbolic
algebra structure is very convenient for the construction of points of 
the spectrum. Many interesting algebras such as the first Weyl 
algebra $A_{1}$, $U(sl_{2})$ and their quantized versions have a 
Hyperbolic algebra structure. Points of the spectrum of the 
category of modules over Hyperbolic algebras are constructed 
in \cite{R}. We recall some basic facts about Hyperbolic algebras 
and two important construction theorems from \cite{R}.

Let $\theta$ be an automorphism of a commutative algebra $R$; and let $\xi$ 
be an element of $R$. 
\begin{defn}
The Hyperbolic algebra $R\{\theta,\,\xi\}$ is defined to be 
the $R-$algebra generated by $x,\,y$ subject to the following relations:
\[
xy=\xi,\quad yx=\theta^{-1}(\xi)
\]
and 
\[
xa=\theta(a)x,\quad ya=\theta^{-1}(a)y
\]
for any $a\in R$. And $R\{\theta,\,\xi\}$ is called a 
Hyperbolic algebra over $R$.
\end{defn}

From \cite{R}, we have the following construction theorems:
\begin{thm}
(Thm 3.2.2.in \cite{R})
\begin{enumerate}
\item Let $P\in Spec(R)$, and assume the orbit of $P$ under the 
action of the automorphism $\theta$ is infinite.
\begin{enumerate}
\item If $\theta^{-1}(\xi)\in P$, and $\xi \in P$, then the left ideal 
\[
P_{1,\,1}\colon=P+R\{\theta,\,\xi\}x+R\{\theta,\,\xi\}y
\]
is a two-sided ideal from $Spec_{l}(R\{\theta,\,\xi\})$.

\item If $\theta^{-1}(\xi)\in P$, $\theta^{i}(\xi)\notin P$ for $0\leq
  i\leq n-1$, 
and $\theta^{n}(\xi)\in P$, then the left ideal 
\[
P_{1,\,n+1}\colon=R\{\theta,\,\xi\}P+R\{\theta,\,\xi\}x+R\{\theta,\,\xi\}y^{n+1}
\]
belongs to $Spec_{l}(R\{\theta,\,\xi\})$.

\item If $\theta^{i}(\xi)\notin P$ for $i\geq 0$ and
  $\theta^{-1}(\xi)\in P$, then
\[
P_{1,\,\infty}\colon=R\{\theta,\,\xi\}P+R\{\theta,\,\xi\}x
\]
belongs to $Spec_{l}(R\{\theta,\,\xi\})$.

\item If $\xi \in P $ and $\theta^{-i}(\xi)\notin P$ for all $i\geq
  1$, 
then the left ideal 
\[
P_{\infty,\,1}\colon=R\{\theta,\,\xi\}P+R\{\theta,\,\xi\}y
\]
belongs to $Spec_{l}(R\{\theta,\,\xi\})$.
\end{enumerate}
\item If the ideal $P$ in (b),\, (c),\, or (d) is maximal, 
then the corresponding left ideal of $Spec_{l}(R\{\theta,\,\xi\})$ is maximal.

\item Every left ideal $Q \in Spec_{l}(R\{\theta,\,\xi\})$ such that
  $\theta^{\nu}(\xi)\in Q$ for 
a $\nu \in \Z$ is equivalent to one left ideal as defined above
  uniquely from a prime ideal $P \in Spec(R)$. The latter means that 
if $P$ and $P'$ are two prime ideals of $R$ and $(\alpha,\,\beta)$ and
  $(\nu,\,\mu)$ 
take values $(1,\,\infty),\,(\infty,\,1),\,(\infty,\,\infty)$ or
  $(1,\,n)$, 
then $P_{\alpha,\,\beta}$ is equivalent to $P'_{\nu,\,\mu}$ if and
  only if $\alpha=\nu,\,\beta=\mu$ and $P=P'$.
\end{enumerate}
\end{thm}
\qed

\begin{thm}
(Prop 3.2.3. in \cite{R}) 
\begin{enumerate}
\item Let $P\in Spec(R)$ be a prime ideal of $R$ 
such that $\theta^{i}(\xi)\notin P$ for $i\in \Z$ and $\theta^{i}(P)-P\neq \O$ 
for $i\neq 0$, then $P_{\infty,\,\infty}=R\{\xi,\,\theta\}P\in Spec_{l}(R\{\xi,\,\theta\})$.
\item Moreover, if ${\bf P}$ is a left ideal of $R\{\theta,\,\xi\}$ such 
that ${\bf P}\cap R=P$, then ${\bf P}=P_{\infty,\,\infty}$. In particular,
if $P$ is a maximal ideal, then $P_{\infty,\,\infty}$ is a maximal 
left ideal.
\item If a prime ideal $P'\subset R$ is such that
  $P_{\infty,\,\infty}=P'_{\infty,\,\infty}$, 
then $P'=\theta^{n}(P)$ for some integer $n$. Conversely,
$\theta^{n}(P)_{\infty,\,\infty}=P_{\infty,\,\infty}$ for any $n\in \Z$. 
\end{enumerate}
\end{thm}\qed

\subsection{Realize $U_{q}(f(K))$ as Hyperbolic algebras}
Let $R$ be the sub-algebra of $U_{q}(f(K))$ generated by 
$EF,\,K^{\pm 1}$, then $R$ is a commutative algebra. We define 
an algebra automorphism of $R$ as follows:
\[
\theta \colon R \longrightarrow R
\]
by 
\[
\theta(EF)=EF+f(\theta(K)), \quad \theta(K^{\pm 1})=q^{\mp 2}K^{\pm 1}.
\]
It is obvious that $\theta$ extends to an algebra automorphism of $R$.
We also have the following lemma:
\begin{lem} The following identities hold:
\begin{gather*}
E(EF)=\theta(EF)E,
\\
F(EF)=\theta^{-1}(EF)F,
\\
EK=\theta(K)E,\\
FK=\theta^{-1}(K)F.
\end{gather*}
\end{lem}
{\bf Proof:} The Proof is straightforward.\qed

From Lemma 2.1, we have the following:
\begin{prop}
$U_{q}(f(K))=R\{\xi=EF,\,\theta\}$ is a Hyperbolic algebra
 with $R$ and $\theta$ defined as above.
\end{prop}
\qed

We have a corollary:
\begin{cor}
The Gelfand-Kirillov dimension of $U_{q}(f(K))$ is $3$.
\end{cor}
{\bf Proof:} This follows from the fact that $R$ has Gelfand-Kirillov
dimension $2$ and $U_{q}(f(K))$ is a Hyperbolic algebra over $R$.\qed

\subsection{Irreducible weight representations of $U_{q}(f_{m}(K))$} 
Now we apply the above construction theorems to the case of
$U_{q}(f_{m}(K))$ and construct families of irreducible weight 
representations for $U_{q}(f_{m}(K))$. 

Let $\alpha, 0\neq \beta \in \C$. We denote by
$P=M_{\alpha,\,\beta}=(\xi-\alpha, \,K-\beta)\subset R$ the maximal
ideal of $R$ generated by $\xi-\alpha, K-\beta$. We need a lemma:
\begin{lem}
$\theta^{n}(M_{\alpha,\,\beta})\neq M_{\alpha,\,\beta}$ for 
any $n \geq 1$. In particular, $M_{\alpha,\,\beta}$ has infinite orbit
under the action of $\theta$.
\end{lem}
{\bf Proof:} We have 
\begin{eqnarray*}
\theta^{n}(K-\beta) & = & (q^{-2n}K-\beta)\\
& =& q^{-2n}(K-q^{2n}\beta).
\end{eqnarray*}
Since $q$ is not a root of unity, $q^{2n}\neq 1 $ for any $n\neq
0$. Thus $\theta^{n}(M_{\alpha,\,\beta})\neq M_{\alpha,\,\beta}$ for any
$n\geq 1$.\qed

Another lemma is in order:
\begin{lem}
\begin{enumerate}
\item For $n\geq 0$, we have the following:
\[
\theta^{n}(EF)=EF+\frac{1}{q-q^{-1}}(\frac{q^{-2m}(1-q^{-2nm})}{1-q^{-2m}}K^{m}-\frac{q^{2m}(1-q^{2nm})}{1-q^{2m}}K^{-m}).
\]
\item For $n\geq 1$, we have the following:
\[
\theta^{-n}(EF)=EF-\frac{1}{q-q^{-1}}(\frac{1-q^{2nm}}{1-q^{2m}}K^{m}-\frac{1-q^{-2nm}}{1-q^{-2m}}K^{-m}).
\]
\end{enumerate}
\end{lem}
{\bf Proof:} For $n\geq 1$, we have 
\begin{eqnarray*}
\theta^{n}(EF)&=&EF+\frac{1}{q-q^{-1}}((q^{-2m}+\cdots+q^{-2nm})K\\
&-&(q^{2m}+\cdots+q^{2nm})K^{-m})\\
&=& EF+\frac{1}{q-q^{-1}}(\frac{q^{-2m}(1-q^{-2nm})}{1-q^{-2m}}K^{m}-\frac{q^{2m}(1-q^{2nm})}{1-q^{2m}}K^{-m}).
\end{eqnarray*}
The second statement can be verified similarly.\qed

\begin{thm}
\begin{enumerate}
\item If $\alpha= \frac{\beta^{m}-\beta^{-m}}{q-q^{-1}},\,
\beta^{m}=\pm q^{mn}$ for some $n\geq 0$, then 
$\theta^{n}(\xi)\in M_{\alpha,\,\beta},\,$ and $\theta^{-1}(\xi) \in
M_{\alpha,\,\beta}$, 
thus $U_{q}(f_{m}(K))/P_{1,\,n+1}$ is a finite dimensional 
irreducible representation of $U_{q}(f(K))$.

\item If $\alpha= \frac{\beta^{m}-\beta^{-m}}{q-q^{-1}}$ and 
$\beta^{m}\neq \pm q^{mn}$ for all $n\geq 0$, 
then $U_{q}(f_{m}(K))/P_{1,\,\infty}$ is an infinite 
dimensional irreducible representation of $U_{q}(f_{m}(K))$.

\item If $\alpha=0 $ and $0\neq \frac{1}{q-q^{-1}}(\frac{1-q^{2nm}}{1-q^{2m}}\beta^{m}-\frac{1-q^{-2nm}}{1-q^{-2m}}\beta^{-m})$ for any $n\geq 1$, 
then $U_{q}(f_{m}(K))/P_{\infty,1}$ is an infinite dimensional 
irreducible representation of $U_{q}(f_{m}(K))$. 
\end{enumerate}
\end{thm}
{\bf Proof:} Since $\theta^{-1}(\xi)=\xi-\frac{K^{m}-K^{-m}}{q-q^{-1}}$, thus $\theta^{-1}(\xi)\in M_{\alpha,\,\beta}$ 
if and only if $\alpha=\frac{\beta^{m}-\beta^{-m}}{q-q^{-1}}$. 
Now by the proof of Lemma 2.3, we have 
\begin{eqnarray*} 
\theta^{n}(\xi)&=&\xi+\frac{1}{q-q^{-1}}((q^{-2m}+\cdots,+q^{-2nm})K^{m}\\
&-&(q^{2m}+\cdots,+q^{2nm})K^{-m})\\
&=& \xi+\frac{1}{q-q^{-1}}(\frac{q^{-2m}(1-q^{-2nm})}{1-q^{-2m}}K^{m}-\frac{q^{2m}(1-q^{2nm})}{1-q^{2m}}K^{-m}).
\end{eqnarray*}
Hence $\theta^{n}(\xi)\in M_{\alpha,\beta}$ if and only if 
\begin{eqnarray*}
0&=&\alpha+\frac{1}{q-q^{-1}}((q^{-2m}+\cdots,+q^{-2nm})\beta^{m}-(q^{2m}+\cdots,+q^{2nm})\beta^{-m})\\
&=& \alpha+\frac{1}{q-q^{-1}}(\frac{q^{-2m}(1-q^{-2nm})}{1-q^{-2m}}\beta^{m}-\frac{q^{2m}(1-q^{2nm})}{1-q^{2m}}\beta^{-m}).
\end{eqnarray*}
Hence when $\alpha= \frac{\beta^{m}-\beta^{-m}}{q-q^{-1}},\,
\beta^{m}=\pm q^{mn}$ for some $n\geq 0$, we have
\[ 
\theta^{n}(\xi)\in M_{\alpha,\,\beta},\, \theta^{-1}(\xi) 
\in M_{\alpha,\,\beta}.
\]
By Theorem 2.1, $U_{q}(f_{m}(K))/P_{1,\,n+1}$ is a finite dimensional 
irreducible representation of $U_{q}(f_{m}(K))$. We have already
proved the first statement, the rest of the statements can be 
similarly verified.\qed

\begin{rem}
The representations we constructed in Theorem 2.3 exhaust all finite
dimensional irreducible representations, highest weight irreducible
representations  and lowest weight irreducible representations of 
$U_{q}(f_{m}(K))$. For each $n\geq 1$, there are exactly $2m$ 
irreducible representations of dimension $n$.
\end{rem}

In addition, we have the following:
\begin{thm}
If $\alpha \neq -\frac{1}{q-q^{-1}}(\frac{q^{-2m}(1-q^{-2nm})}{1-q^{-2m}}\beta^{m}-\frac{q^{2m}(1-q^{2nm})}{1-q^{2m}}\beta^{-m})$ for any $n\geq 0$ and
 $\alpha \neq
 \frac{1}{q-q^{-1}}(\frac{1-q^{-2nm}}{1-q^{-2m}}\beta^{m}-\frac{1-q^{2nm}}{1-q^{2m}}\beta^{-m})$ for any $n\geq 1$, then
 $U_{q}(f_{m}(K))/P_{\infty,\,\infty}$ is an infinite dimensional 
irreducible weight representation of $U_{q}(f_{m}(K))$.
\end{thm}
{\bf Proof:} The proof is similar to the one of Theorem 2.3, we omit it here.\qed
\begin{cor}
The representations constructed in Theorem 2.3 and Theorem 2.4 exhaust all
irreducible weight representations of $U_{q}(f_{m}(K))$.
\end{cor}
{\bf Proof:} It follows directly from Theorems 2.1, 2.2, 2.3 and 2.4.\qed
\begin{rem}
When $m=1$, the above results recover the weight representations of 
$U_{q}(sl_{2})$. So our results are just a natural generalization of 
those for $U_{q}(sl_{2})$.
\end{rem}
\section{The Whittaker model of the center $Z(U_{q}(f_{m}(K)))$}
Let $\frak g$ be a finite dimensional complex semisimple Lie algebra 
and $U(\frak g)$ be its universal enveloping algebra. The Whittaker 
model for the center of $U(\frak g)$ was studied by Kostant in
\cite{K}. The Whittaker model for the center $Z(U(\frak g))$ is 
defined by a non-singular character of the nilpotent Lie subalgebra 
$\frak n^{+}$ of $\frak g$. Using the Whittaker model, Kostant studied 
the structure of Whittaker modules of $U(\frak g)$ and several important 
results about Whittaker modules were obtained in \cite{K}. 
Later on, Kostant's idea was further generalized by Lynch in
\cite{L} and by Macdowell in \cite{M} to the case of singular 
characters of $\frak n^{+}$ and similar results were proved to 
hold. 

The obstacle of generalizing the Whittaker model to the quantized 
enveloping algebra $U_{q}(\frak g)$ with $\frak g$ of higher ranks 
is that there is no non-singular character existing for the positive 
part $(U_{q}(\frak g))^{> 0}$ of $U_{q}(\frak g)$ because of the 
quantum Serre relations. In order to overcome this difficulty, it was 
Sevostyanov who first realized to use the topological version $U_{h}(\frak g)$ 
over $C[[h]]$ of quantum groups. Using a family of Coxeter realizations 
$U^{s_{\pi}}_{h}(\frak g)$ of the quantum group $U_{h}(\frak g)$
indexed by the Coxeter elements $s_{\pi}$, he was able to prove 
Kostant's results for $U_{h}(\frak g)$ in \cite{S1}. While in the 
simplest case of ${\frak g}=sl_{2}$, the quantum Serre relations 
are trivial, and a direct generalization has been found in \cite{O}.

In addition, it is reasonable to expect that the Whittaker model 
exists for most of the deformations of $U_{q}(sl_{2})$. In this 
section, we show that there is such a Whittaker model for the 
center of $U_{q}(f_{m}(K))$ and study their Whittaker modules. 
We prove several similar results as in \cite{K} and \cite{O}. 
We follow the approach in \cite{K} and \cite{O}. For the reader's 
convenience, we will work out all the details here. For the simplicity
of notations, we denote $f_{m}(K)=\frac{K^{m}-K^{-m}}{q-q^{-1}}$ 
by $f(K)$ from now on.

\subsection{The center $Z(U_{q}(f(K)))$ of $U_{q}(f(K))$}
We first give a description of the center of $U_{q}(f(K))$. To
proceed, we define a Casimir element $\Omega$ by: 
\[ 
\Omega=FE+\frac{q^{2m}K^{m}+K^{-m}}{(q^{2m}-1)(q-q^{-1})}.
\]

We have the following:
\begin{prop}
\begin{eqnarray*}
\Omega &=& FE+\frac{q^{2m}K^{m}+K^{-m}}{(q^{2m}-1)(q-q^{-1})}\\
  &=& EF+\frac{K^{m}+q^{2m}K^{-m}}{(q^{2m}-1)(q-q^{-1})}.
\end{eqnarray*}
\end{prop}
{\bf Proof:} The proof is easy.\qed

In addition, we have the following lemma:
\begin{lem}
$\Omega$ is in the center of $U_{q}(f(K))$.
\end{lem}
{\bf Proof:} The proof is easy.\qed

We have the following description of the center $Z(U_{q}(f(K)))$ of 
$U_{q}(f(K))$:
\begin{prop}
(See also \cite{JWZ}) $Z(U_{q}(f(K)))$ is the subalgebra of $U_{q}(f(K))$ 
generated by $\Omega$. In particular, $Z(U_{q}(f(K)))$ is 
isomorphic to the polynomial algebra in one variable.
\end{prop}
{\bf Proof:} 
By Lemma 3.1., we have $\Omega \in Z(U_{q}(f(K)))$. Thus 
the subalgebra $\C[\Omega]$ generated by $\Omega$ is contained 
in $Z(U_{q}(f(K)))$. So it suffices to prove 
the other inclusion $Z(U_{q}(f(K)))\subseteq \mathbb{C}[\Omega]$. 
Note that $U_{q}(f(K))=\bigoplus_{n\in \Z}U_{q}(f(K))_{n}$ 
where $U_{q}(f(K))_{n}$ is the $\C-$span of elements
$\{u\in U_{q}(f(K))\mid Ku=q^{2n}uK\}$. Suppose $x\in Z(U_{q}(f(K)))$, 
then $xK=Kx$. Thus $x\in U_{q}(f(K))_{0}$ which is generated by $EF,
K^{\pm 1}$. By the definition of $\Omega$, we know that
$U_{q}(f(K))_{0}$ is also generated by $\Omega,K^{\pm 1}$. 
Hence $x=\sum f_{i}(\Omega)K^{i}$ where $f_{i}(\Omega)$ are 
polynomials in $\Omega$. Therefore 
\[
xE=\sum f_{i}(\Omega)K^{i}E=\sum f_{i}(\Omega)q^{2i}EK^{i}=Ex
\]
which forces $i=0$. So $x\in \mathbb{C}[\Omega]$ as desired. 
So we have proved that $Z(U_{q}(f(K)))=\mathbb{C}[\Omega]$.
\qed

\subsection{The Whittaker Model for $Z(U_{q}(f(K)))$}
Now we construct the Whittaker model for $Z(U_{q}(f(K)))$ following 
\cite{K} and \cite{O}. For the rest of this paper, we use the term of 
modules instead of representations.

We denote by $U_{q}(E)$ the subalgebra of $U_{q}(f(K))$ generated by
$E$, by $U_{q}(F,\,K^{\pm 1})$ the subalgebra of $U_{q}(f(K))$
generated by $F,\,K^{\pm 1}$.
\begin{defn}
An algebra homomorphism $\eta\colon U_{q}(E)\longrightarrow \mathbb{C}$ 
is called a non-singular character of $U_{q}(E)$ if $\eta(E)\neq 0$.
\end{defn}

From now on, we always fix such a non-singular character of 
$U_{q}(E)$ and denote it by $\eta$. As in \cite{K}, we have 
the following:
\begin{defn}
Let $V$ be a $U_{q}(f(K))-$module, a vector $0\neq v\in V$ is 
called a Whittaker vector of type $\eta$ if $E$ acts on $v$ through 
the non-singular character $\eta$, i.e., $Ev=\eta(E)v$. If
$V=U_{q}(f(K))v$, then we call $V$ a Whittaker module of type $\eta$ 
and $v$ is called a cyclic Whittaker vector of type $\eta$ for $V$. 
\end{defn}

By the definition of $U_{q}(f(K))$, the following decomposition 
of $U_{q}(f(K))$ is obvious:
\begin{prop}
$U_{q}(f(K))$ is isomorphic to $U_{q}(F,K^{\pm 1})\otimes_{\mathbb{C}} 
U_{q}(E)$ as vector spaces and $U_{q}(f(K))$ is a free module over 
the subalgebra $U_{q}(E)$.
\end{prop}
Let us denote the kernel of $\eta \colon U_{q}(E)\longrightarrow
\mathbb{C}$ by $U_{q,\,\eta}(E)$, and we have the following 
decompositions of $U_{q}(E)$ and $U_{q}(f(K))$.
\begin{prop} We have $U_{q}(E)=\mathbb{C} \oplus U_{q,\,\eta}(E)$. 
Furthermore, 
\[
U_{q}(f(K))
\cong U_{q}(F,\,K^{\pm 1})\oplus U_{q}(f(K))U_{q,\,\eta}(E).
\]
\end{prop}
{\bf Proof:} It is obvious that $U_{q}(E)=\mathbb{C} \oplus U_{q,\,\eta}(E)$. 
And we have 
\[
U_{q}(f(K))=U_{q}(F,\,K^{\pm 1})\otimes(\C \oplus U_{q,\,\eta}(E)),
\]

thus 
\[
U_{q}(f(K))\cong U_{q}(F,\,K^{\pm 1})\oplus U_{q}(f(K))U_{q,\,\eta}(E)
\]
So we are done.
\qed

Now we define a projection 
\[
\pi\colon U_{q}(f(K))\longrightarrow U_{q}(F,\,K^{\pm 1})
\]
 from $U_{q}(f(K))$ onto $U_{q}(F,\,K^{\pm 1})$ by taking 
the $U_{q}(F,\,K^{\pm 1})-$component of any $u\in U_{q}(f(K))$. 
We denote the image $\pi(u)$ of $u\in U_{q}(f(K))$ by $u^{\eta}$ for short.
\begin{lem}
If $v \in Z(U_{q}(f(K)))$, then we have $u^{\eta}v^{\eta}=(uv)^{\eta}$
for any $u\in U_{q}(f(K))$.
\end{lem}
{\bf Proof:} Let $u\in U_{q}(f(K)),\,v \in Z(U_{q}(f(K)))$, then we have
\begin{eqnarray*}
 uv-u^{\eta}v^{\eta}&=&(u-u^{\eta})v+u^{\eta}(v-v^{\eta})\\
&=& v(u-u^{\eta})+u^{\eta}(v-v^{\eta}),
\end{eqnarray*}
which is in $U_{q}(f(K))U_{q,\,\eta}(E)$. Hence $(uv)^{\eta}=u^{\eta}v^{\eta}$.
\qed

From the definition of $\Omega$, we have the following description of $\pi(\Omega)$:
\begin{lem}
\[
\pi(\Omega)=\eta(E)F+\frac{q^{2m}K^{m}+K^{-m}}{(q^{2m}-1)(q-q^{-1})}.
\]
\end{lem}\qed

\begin{prop}
The map 
\[
\pi\colon Z(U_{q}(f(K))\longrightarrow U_{q}(F,\,K^{\pm 1})
\]
is an algebra isomorphism of $Z(U_{q}(f(K)))$ onto 
its image $W(F,\,K^{\pm 1})$ in $U_{q}(F,\,K^{\pm 1})$.
\end{prop}
{\bf Proof:} It follows from Lemma 3.2. that $\pi$ is a
homomorphism of algebras. It remains to show that $\pi$ 
is injective. Suppose $\pi(u)=0$ for some $0\neq u\in Z(U_{q}(f(K)))$.
Since $Z(U_{q}(f(K)))$ is a polynomial algebra in $\Omega$, then
$u=\sum_{i=0}^{N}u_{i}\Omega^{i}$. By Lemma 3.4., 
$\pi(\Omega)=\eta(E)F+\frac{q^{2m}K^{m}+K^{-m}}{(q^{2m}-1)(q-q^{-1})}$.
Thus
\[
\pi(u)=\sum_{i=0}^{N}u_{i}(\eta(E)F+\frac{q^{2m}K^{m}+K^{-m}}{(q^{2m}-1)(q-q^{-1})})^{i}=0.
\]
By direct computations, we have $u_{N}(\eta(E))^{N}F^{N}=0$. So $u_{N}(\eta(E))^{N}=0$, which is a contradiction. So $\pi$ is an injection. Thus 
$\pi$ is an algebra isomorphism from $Z(U_{q}(f(K)))$ onto 
its image $W(F,\,K^{\pm 1})$ in $U_{q}(F,\,K^{\pm 1})$.\qed

\begin{lem}
If $u^{\eta}=u$, then we have 
\[
u^{\eta}v^{\eta}=(uv)^{\eta}
\]
 for any $v\in U_{q}(f(K))$.
\end{lem}
{\bf Proof:} We have 
\begin{eqnarray*}
uv-u^{\eta}v^{\eta}&=&(u-u^{\eta})v+u^{\eta}(v-v^{\eta})\\
&=&u^{\eta}(v-v^{\eta})
\end{eqnarray*} 
which is in $U_{q}(f(K))U_{q,\eta}(E)$. So we have 
\[
u^{\eta}v^{\eta}=(uv)^{\eta}
\]
for any $v\in U_{q}(f(K))$.
\qed

Let $\tilde{A}$ be the subalgebra of $U_{q}(f(K))$ 
generated by $K^{\pm 1}$. Then $\tilde{A}$ is a graded 
vector space with $A_{0}=\C$, and 
\[
\tilde{A}_{[n]}=\C K^{n}\oplus \C K^{-n}
\]
for $n\geq 1$, and 
\[
\tilde{A}_{[n]}=0
\]
for $n\leq -1$. 

As in \cite{K} and \cite{O}, we define a filtration of $U_{q}(F,\,K^{\pm
1})$ 
as follows:
\[
U_{q}(F,\,K^{\pm 1})_{[n]}=\bigoplus_{im+\mid j \mid \leq nm} 
U_{q}(F,\,K^{\pm 1})_{i,j}
\]
where $U_{q}(F,\,K^{\pm 1})_{i,j}$ is the vector space 
spanned by $F^{i}K^{j}$. We denote by 
\[
W(F,\,K^{\pm 1})_{[p]}=\C-span\{1,\,\Omega^{\eta},\,\cdots,\,(\Omega^{\eta})^{p}\}
\]
for $p\geq 0$. It is easy to see that
\[
W(F,\,K^{\pm 1})_{[p]}\subset W(F,\,K^{\pm 1})_{[p+1]},
\quad W(F,\,K^{\pm 1})=\sum_{p\geq 0}W(F,\,K^{\pm 1})_{[p]}
\]
And $W(F,\,K^{\pm 1})_{[p]}$ give a filtration 
of $W(F,\,K^{\pm 1})$, which is compatible with 
the filtration of $U_{q}(F,\,K^{\pm 1})$. In particular,
\[
W(F,\,K^{\pm 1})_{[p]}=W(F,\,K^{\pm 1})\cap U_{q}(F,\,K^{\pm 1})_{[p]}
\]
for $q\geq 0$.

Now we have the following decomposition 
of $U_{q}(F,\,K^{\pm 1})$.
\begin{thm}
$U_{q}(F,\,K^{\pm 1})$ is free (as a right module) 
over $W(F,\,K^{\pm 1})$. And the multiplication induces an isomorphism 
\[
\Phi\colon \tilde{A}\otimes W(F,\,K^{\pm 1})\longrightarrow
U_{q}(F,\,K^{\pm 1})
\]
as right $W(F,\,K^{\pm 1})-$modules. In particular, we have the following
\[
\bigoplus_{l+pm=nm}\tilde{A}_{[l]}\otimes W(F,\,K^{\pm 1})_{[p]}\cong
U_{q}(F,\,K^{\pm 1})_{[n]} 
\]
\end{thm}
{\bf Proof:} First of all, the map 
$\tilde{A}\times W(F,\,K^{\pm 1}) \longrightarrow U_{q}(F,\,K^{\pm 1})$ 
is bilinear. So by the universal property of 
the tensor product, there is a map 
from $\tilde{A}\otimes W(F,\,K^{\pm 1})$ 
into $U_{q}(F,\,K^{\pm 1})$ defined by the 
multiplication. It is easy to check 
this map is a homomorphism of right 
$W(F,\,K^{\pm 1})-$modules and is surjective 
as well. Now we show that it is injective.
Let $0\neq u\in \tilde{A}\otimes W(F,\,K^{\pm 1})$ such that $\Phi(u)=0$. We write
$u=\sum_{i=0}^{L}a_{i}(K,\,K^{-1})\otimes b_{i}(\pi(\Omega))^{i}$, where
$a_{i}(K,\,K^{-1})$ are non-zero Laurent polynomials in $\C[K,\,K^{-1}]$ and
$b_{i}\in \C^{\ast}$. Then
$0=\Phi(u)=\sum_{i=0}^{L}a_{i}(K,\,K^{-1})b_{i}(\eta(E)F+\frac{q^{2m}K^{m}+K^{-m}}{(q^{2m}-1)(q-q^{-1})})^{i}$.
By direct computations, we have so $a_{L}(K,\,K^{-1})b_{L}(\eta(E))^{L}F^{L}=0$. 
Thus $a_{L}(K,\,K^{-1})b_{L}(\eta(E))^{L}=0$, which is a contradiction. So 
we have proved that $\Phi$ is a isomorphism of vector spaces. 
In addition, by counting the degrees of both sides, we also have 
\[
\bigoplus_{l+pm=nm}\tilde{A}_{[l]}\otimes W(F,\,K^{\pm 1})_{[p]}\cong
U_{q}(F,\,K^{\pm 1})_{[n]} 
\]
Thus we have proved the theorem.
\qed

Let $Y_{\eta}$ be the left $U_{q}(f(K))-$module 
defined by 
\[
Y_{\eta}=U_{q}(f(K))\otimes_{U_{q}(E)} \C_{\eta}
\]
where $\C_{\eta}$ is one dimensional $U_{q}(E)-$module 
defined by the character $\eta$. It is easy to see that 
\[
Y_{\eta}\cong U_{q}(f(K))/U_{q}(f(K))U_{q,\,\eta}(E)
\]
is a Whittaker module with a cyclic vector denoted by $1_{\eta}$. 
Now we have a quotient map from $U_{q}(f(K))$ to $Y_{\eta}$ 
\[
U_{q}(f(K))\longrightarrow Y_{\eta}.
\]
If $u\in U_{q}(f(K))$, then there is a $u^{\eta}$ which is 
the unique element in $U_{q}(F,\,K^{\pm 1})$ such that 
$u 1_{\eta}=u^{\eta}1_{\eta}$. As in \cite{K}, we define 
the $\eta -$reduced action of $U_{q}(E)$ on $U(F,\,K^{\pm 1})$ 
as follows:
\[
x\bullet v=(xv)^{\eta}-\eta(x)v
\]
where $x\in U_{q}(E)$ and $v\in U_{q}(F,\,K^{\pm 1})$.
\qed
\begin{lem}
Let $u\in U_{q}(F,\,K^{\pm 1})$ and $x\in U_{q}(E)$, 
we have
\[
x\bullet u^{\eta}=[x,\,u]^{\eta}
\]
\end{lem}
{\bf Proof:} $[x,\,u]1_{\eta}=(xu-ux)1_{\eta}=(xu-\eta(x)u)1_{\eta}$. Hence
\[[x,\,u]^{\eta}=(xu)^{\eta}-\eta(x)u^{\eta}=(xu^{\eta})^{\eta}-\eta(x)u^{\eta}=x\bullet u^{\eta}\].\qed

\begin{lem}
Let $x\in U_{q}(E)$, $u\in U_{q}(F,\,K^{\pm 1})$, and $v\in
W(E,\,K^{\pm 1})$, then we have 
\[
x \bullet (uv)=(x\bullet u)v.
\]
\end{lem}
{\bf Proof:} Let $v=w^{\eta}$ for some $w\in Z(U_{q}(f(K))$, 
then $uv=uw^{\eta}=u^{\eta}w^{\eta}=(uw)^{\eta}$. 
Thus
\begin{eqnarray*}
x\bullet(uv)&= &x\bullet(uw)^{\eta}=[x,\,uw]^{\eta}\\
&=&([x,\,u]w)^{\eta}=[x,\,u]^{\eta}w^{\eta}\\
&=&[x,\,u]^{\eta}v\\
&=&(x\bullet u^{\eta})v\\
&=&(x\bullet u)v.
\end{eqnarray*}
So we are done.\qed

Let $V$ be an $U_{q}(f(K))-$module and let $U_{q,\,V}(f(K))$ be the
annihilator of $V$ in $U_{q}(f(K))$. Then $U_{q,\,V}(f(K))$ defines a central 
ideal $Z_{V}\subset Z(U_{q}(f(K)))$ by setting $Z_{V}=U_{q,\,V}(f(K))\cap Z(U_{q}(f(K)))$.
Suppose that $V$ is a Whittaker module with a cyclic Whittaker vector $w$, 
we denote by $U_{q,\,w}(f(K))$ the annihilator of $w$ in $U_{q}(f(K))$. 
It is obvious that 
\[
U_{q}(f(K))U_{q,\,\eta}(E)+U_{q}(f(K))Z_{V}\subset U_{q,\,w}(f(K)).
\]

In the next theorem we show that the reverse inclusion holds.
First of all, we need an auxiliary Lemma:
\begin{lem}
Let $X=\{v\in U_{q}(F,\,K^{\pm 1})\mid (x \bullet v)w=0,\,x\in U_{q}(E)\}$.
Then
\[
X=\tilde{A}\otimes W_{V}(F,\,K^{\pm 1})+W(F,\,K^{\pm 1})
\]
where $W_{V}(F,\,K^{\pm 1})=(Z_{V})^{\eta}$. In fact, 
$U_{q,V}(F,\,K^{\pm 1})\subset X$ and 
\[
U_{q,w}(F,\,K^{\pm 1})=\tilde{A}\otimes W_{w}(F,\,K^{\pm 1})
\]
where $U_{q,\,w}(F,\,K^{\pm 1})=U_{q,\,w}(f(K))\cap U_{q}(F,\,K^{\pm 1})$
\end{lem}
{\bf Proof:} Let us denote by $Y=\tilde{A}\otimes W_{V}(F,\,K^{\pm
  1})+W(F,\,K^{\pm 1})$ where $W(F,\,K^{\pm 1})=(Z(U_{q}(f(K))))^{\eta}$. 
Thus we need to verify $X=Y$. Let $v\in W(F,\,K^{\pm 1})$, 
then $v=u^{\eta}$ for some $u\in Z(U_{q}(f(K)))$. 
So we have 
\begin{eqnarray*}
x\bullet v &=&x\bullet u^{\eta}\\
&=&[x,u]^{\eta}\\
&=&(xu)^{\eta}-\eta(x)u^{\eta}\\
&=&x^{\eta}u^{\eta}-\eta(x)u^{\eta}\\
&=&0.
\end{eqnarray*}
So we have $W(F,\,K^{\pm 1})\subset X$.
Let $u\in Z_{V}$ and $v\in U_{q}(F,\,K^{\pm 1})$. 
Then for any $x\in U_{q}(F)$, we have 
\[
x\bullet(vu^{\eta})=(x\bullet v)u^{\eta}
\]
Since $u\in Z_{V}$, then $u^{\eta}\in U_{q,\,w}(f(K))$. 
Thus we have $vu^{\eta}\in X$, hence 
\[
\tilde{A}\otimes W_{V}(F,\,K^{\pm 1})\subset X
\]
which proves $Y\subset X$. Note that $\tilde{A_{[i]}}$ 
is the two dimensional subspace of $\C[K^{\pm 1}]$ spanned by $K^{\pm i}$ 
and $\overline{W_{V}(F,\,K^{\pm 1})}$ is the complement of 
$W_{V}(F,\,K^{\pm 1})$ in $W(F,\,K^{\pm 1})$. 
Let us set 
\[
M_{i}=\tilde{A_{[i]}}\otimes \overline{W_{V}(F,\,K^{\pm 1})}
\]
thus we have the following:
\[
U_{q}(F,\,K^{\pm 1})=M\oplus Y
\]
where $M=\sum_{i\geq 1}M_{i}$.
We show that $M\cap X\neq 0$. 
Let $M_{[k]}=\sum_{1\leq i\leq k}M_{i}$, 
then $M_{[k]}$ are a filtration of $M$. 
Suppose $n$ is the smallest $n$ such that $X\cap M_{[n]}\neq 0$ 
and $0\neq y\in X\cap M_{[n]}$. 
Then we have $y=\sum_{1\leq i\leq n}y_{i}$ 
where $y_{i}\in \tilde{A_{i}}\otimes \overline{W_{V}(F,\,K^{\pm 1})}$. 
Suppose we have chosen $y$ in such a way that $y$ has the fewest
terms. Then by direct computations, we have $0\neq
y-\frac{1}{\eta(E)(q^{-2n}-1)} E \bullet y\in X\cap M_{[n]}$ with fewer
terms than $y$. This is a contradiction. So we have $X\cap M=0$.

Now we prove that $U_{q,\,w}(F,\,K^{\pm 1})\subset X$. 
Let $u\in U_{q,\,w}(F,\,K^{\pm q})$ and $x\in U_{q}(E)$, 
then we have $xuw=0$ and $uxw=\eta(x)uw=0$. 
Thus $[x,\,u]\in U_{q,\,w}(F,\,K^{\pm 1})$, 
hence $[x,\,u]^{\eta}\in U_{q,\,w}(F,\,K^{\pm 1})$. 
Since $u\in U_{q,\,w}(F,\,K^{\pm 1})\subset U_{q,\,w}(E,\,F,\,K^{\pm 1})$, 
then $x\bullet u=[x,\,u]^{\eta}$. 
Thus $x\bullet u \in U_{q,w}(F,\,K^{\pm 1})$. 
So $u\in X$ by the definition of $X$. 
Now we are going to prove the following:
\[
W(F,\,K^{\pm 1})\cap U_{q,\,w}(F,\,K^{\pm 1})=W_{V}(F,\,K^{\pm 1})
\] 
In fact, $W_{V}(F,\,K^{\pm 1})=(Z_{V}^{\eta})$ and 
$W_{V}(F,\,K^{\pm 1})\subset U_{q,\,w}(F,\,K^{\pm 1})$. 
So if $v\in W_{w}(F,\,K^{\pm 1})\cap U_{q,w}(F,\,K^{\pm 1})$, 
then we can uniquely write $v=u^{\eta}$ for $u\Z(U_{q}(f(K)))$. 
Then $vw=0$ implies $uw=0$ and hence 
$u\in Z(U_{q}(f(K)))\cap U_{q,w}(F,\,K^{\pm 1})$. 
Since $V$ is cyclically generated by $w$, we 
have proved the above statement.
Obviously, we have 
$U_{q}(f(K))Z_{V}\subset U_{q,\,w}(f(K))$. 
Thus we have 
$\tilde{A}\otimes W_{V}(F,\,K^{\pm 1})\subset U_{q,\,w}(F,\,K^{\pm 1})$, 
hence we have $U_{q,w}(F,\,K^{\pm 1})=\tilde{A}\otimes W_{V}(F,\,K^{\pm 1})$.
So we have finished the proof.\qed
\begin{thm}
Let $V$ be a Whittaker module admitting a cyclic 
Whittaker vector $w$, then we have 
\[
U_{q,\,w}(f(K))=U_{q}(f(K))Z_{V}+U_{q}(f(K))U_{q,\,\eta}(E).
\]
\end{thm}
{\bf Proof:} It is obvious that 
\[
U_{q}(f(K))Z_{V}+U_{q}(f(K))U_{\eta}(E)
\subset U_{q,\,w}(f(K))
\]
Let $u\in U_{q,\,w}(f(K))$, we show that 
$u\in U_{q}(f(K))Z_{V}+U_{q}(f(K))U_{q,\,\eta}(E)$. 
Let $v=u^{\eta}$, then it suffices to show 
that $v\in \tilde{A}\otimes W_{V}(F,\,K^{\pm 1})$. 
But $v\in U_{q,\,w}(F,\,K^{\pm 1})=\tilde{A}\otimes W_{V}(F,\,K^{\pm 1})$.
So we have proved the theorem.\qed
\begin{thm}
Let $V$ be any Whittaker module for $U_{q}(f(K))$, 
then the correspondence
\[
V \longrightarrow Z_{V}
\]
sets up a bijection between the set of all 
equivalence classes of Whittaker modules and 
the set of all ideals of $Z(U_{q}(f(K)))$.
\end{thm}
{\bf Proof:} Let $V_{i}, i=1,2$ be two Whittaker modules. 
If $Z_{V_{1}}=Z_{V_{2}}$, then clearly $V_{1}$ is 
equivalent to $V_{2}$ by the above Theorem. 
Now let $Z_{\ast}$ be an ideal of $Z(U_{q}(f(K)))$ 
and let $L=U_{q}(f(K))Z_{\ast}+U_{q}(f(K))U_{\eta}(E)$. 
Then $V=U_{q}(f(K))/L$ is a Whittaker module 
with a cyclic Whittaker vector $w=\bar{1}$. 
Obviously we have $U_{q,w}(f(K))=L$. So that 
$L=U_{q,\,w}(f(K))=U_{q}(f(K))Z_{V}+U_{q}(f(K))U_{q,\,\eta}(E)$. 
This implies that 
\[
\eta(L)=\pi(Z_{\ast})=\pi(R_{w})=\pi(Z_{V}).
\]
Since $\pi$ is injective, thus $Z_{V}=Z_{\ast}$. Thus we have finished 
the proof.\qed

\begin{thm}
Let $V$ be an $U_{q}(f(K))-$module. Then $V$ is a Whittaker module if 
and only if 
\[
V\cong U_{q}(f(K))\otimes_{Z(U_{q}(f(K)))
\otimes U_{q}(E)}(Z(U_{q}(f(K)))/Z_{\ast})_{\eta}.
\]
In particular, in such a case the ideal $Z_{\ast}$ 
is uniquely determined to be $Z_{V}$.
\end{thm}
{\bf Proof:} If $1_{\ast}$ is the image 
of $1$ in $Z(U_{q}(f(K)))/Z_{\ast}$, 
then 
\[
Ann_{Z(U_{q}(f(K)))\otimes U_{q}(F)}(1_{\ast})
=U_{q}(E)Z_{\ast}+Z(U_{q}(f(K)))U_{q,\eta}(E).
\]
Thus the annihilator of $w=1\otimes 1_{\ast}$ 
is 
\[
U_{q,\,w}(f(K))=U_{q}(f(K))Z_{\ast}+U_{q}(f(K))U_{q,\,\eta}(E).
\] 
Then the result follows from the previous theorem. \qed

\begin{thm}
Let $V$ be an $U_{q}(f(K))-$module with a cyclic Whittaker vector
$w\in V$. Then any $v\in V$ is a Whittaker vector if and only if
$v=uw$ for some $u\in Z(U_{q}(f(K)))$.
\end{thm}
{\bf Proof:} If $v=uw$ for some $u \in Z(U_{q}(f(K)))$, 
then it is obvious that $v$ is a Whittaker vector. 
Conversely, let $v=uw$ for some $u\in U_{q}(f(K))$ 
be a Whittaker vector of $V$. Then $v=u^{\eta}w$ by 
the definition of Whittaker module. So we may 
assume that $u\in U_{q}(F,\,K^{\pm 1})$. If $x\in U_{q}(E)$, 
then we have $xuw=\eta(x)uw$ and $uxw=\eta(x)uw$. 
Thus $[x,\,u]w=0$ and hence $[x,\,u]^{\eta}w=0$. But 
we have $x\bullet u=[x,\,u]^{\eta}$. Thus we have 
$u\in X$. We can now write $u=u_{1}+u_{2}$ with 
$u_{1}\in U_{q}(F,\,K^{\pm 1})$ and $u_{2}\in W(F,\,K^{\pm 1})$. 
Then $u_{1}w=0$. Hence $u_{2}w=v$. But $u_{2}=u_{3}^{\eta}$ 
with $u_{3}\in Z(U_{q}(f(K)))$. So we have $v=u_{3}w$, which 
proves the theorem.
\qed
 
Now let $V$ be a Whittaker module and $End_{U_{q}(f(K))}(V)$ be the 
endomorphism ring of $V$ as a $U_{q}(f(K))-$module. Then we can define 
the following homomorphism of algebras using the action of $Z(U_{q}(f(K)))$ 
on $V$:
\[
\pi_{V}\colon Z(U_{q}(f(K))\longrightarrow End_{U_{q}(f(K))}(V)
\]
It is clear that 
\[
Z(U_{q}(f(K)))/Z_{V}(U_{q}(f(K))) \cong \pi_{V}(Z(U_{q}(f(K)))) \subset End_{U_{q}(f(K))}(V).
\]
In fact, the next theorem says that 
this inclusion is an equality as well.
\begin{thm}
Let $V$ be a Whittaker $U_{q}(f(K))-$module. Then 
$End_{U_{q}(f(K))}(V)\cong Z(U_{q}(f(K)))/Z_{V}$. 
In particular, $End_{U_{q}(f(K))}(V)$ is commutative.
\end{thm}
{\bf Proof:} Let $w\in V$ be a cyclic Whittaker vector. 
If $\alpha \in End_{U_{q}(f(K))}(V)$, then 
$\alpha(w)=uw$ for some $u\in Z(U_{q}(f(K)))$ by Theorem 3.5. 
Thus we have $\alpha(vw)=vuw=uvw$. Hence $\alpha=\pi_{u}$, which 
proves the theorem.
\qed

We have the following description about the basis of an 
irreducible Whittaker module $(V,\,w)$ where $w\in V$ is a cyclic 
Whittaker vector. 
\begin{thm}
Let $(V,\,w)$ be an irreducible Whittaker module 
with a Whittaker vector $w$, then $V$ has a 
$\C-$basis consisting of elements $\{K^{i}w \mid i\in \Z \}$.
\end{thm}
{\bf Proof:} The proof is straightforward.
\qed

Now we are going to construct explicitly some Whittaker 
modules. Let 
\[
\xi \colon Z(U_{q}(f(K))) \longrightarrow \C
\] 
be a central character of the center $Z(U_{q}(f(K)))$. For any given 
central character $\xi$, let $Z_{\xi}=Ker(\xi)\subset Z(U_{q}(f(K)))$ 
and $Z_{\xi}$ is a maximal ideal of $Z(U_{q}(f(K)))$. Since $\C$ is 
algebraically closed, then $Z_{\xi}=(\Omega -a_{\xi})$ for some 
$a_{\xi} \in \C$. For any given central character $\xi$, let $\C_{\xi,\eta}$ 
be the one dimensional $Z(U_{q}(f(K)))\otimes U_{q}(E)-$module defined 
by $uvy=\xi(u)\eta(v)y$ for any $u\in Z(U_{q}(f(K)))$ and any $v\in
U_{q}(E)$. We set
\[
Y_{\xi,\,\eta}=U_{q}(f(K))\otimes_{Z(U_{q}(f(K)))\otimes U_{q}(E)} \C_{\xi,\,\eta}.
\] 
It is easy to see that $Y_{\xi,\,\eta}$ is a Whittaker module of type
$\eta$ and admits a central character $\xi$. By Schur's lemma, we know 
every irreducible representation has a central character. As studied
in \cite{JWZ}, we know $U_{q}(f(K))$ has a similar theory for Verma
modules. In fact, Verma modules also fall into the category of
Whittaker modules if we take the trivial character of $U_{q}(E)$. 
Namely we have the following 
\[
M_{\lambda}=U_{q}(f(K))\otimes_{U_{q}(E,\,K^{\pm 1})}\C_{\lambda}
\]
where $K$ acts on $\C_{\lambda}$ through the character $\lambda$ of
$\C[K^{\pm 1}]$ and $U_{q}(E)$ acts trivially on $\C_{\lambda}$. It is
obvious that $M_{\lambda}$ admits a central character. It is
well-known that Verma modules may not be necessarily irreducible, 
even though they have central characters. While Whittaker modules 
are in the other extreme as shown in the next theorem:
\begin{thm}
Let $V$ be a Whittaker module for $U_{q}(f(K))$. Then the following
statements are equivalent.\\
\begin{enumerate}
\item $V$ is irreducible.\\
\item $V$ admits a central character.\\
\item $Z_{V}$ is a maximal ideal.\\
\item The space of Whittaker vectors of $V$ is one-dimensional.\\
\item All nonzero Whittaker vectors of $V$ are cyclic.\\
\item The centralizer $End_{U_{q}(f(K))}(V)$ is reduced to  $\C$.\\
\item $V$ is isomorphic to $Y_{\xi,\eta}$ for some central character $\xi$.
\end{enumerate}
\end{thm}
{\bf Proof:} It is easy to see that $(2)-(7)$ are equivalent 
to each other by using the previous Theorems we have just proved. 
Since $\C$ is algebraically closed and uncountable, we also know $(1)$ 
implies $(2)$ by using a theorem due to Dixmier \cite{Di}. To complete 
the proof, it suffices to show that $(2)$ implies $(1)$, namely if $V$ has
a central character, then $V$ is irreducible. Let $\omega \in V$ be a
cyclic Whittaker vector, then $V=U_{q}(f(K))\omega$. Since $V$ has a 
central character, then it is easy to see from the description of the
center that $V$ has a $\C-$basis consisting of elements
$\{K^{i}\omega \mid i \in \Z\}$. Let $0\neq
v=\sum_{0}^{n}a_{i}K^{i}\omega \in V$, then
$E\sum_{0}^{n}a_{i}K^{i}\omega=\sum_{0}^{n} q^{-2i}a_{i}K^{i}E\omega=\eta(E)\sum_{0}^{n}q^{-2i}a_{i}K^{i}\omega$. Thus we have $0 \neq \eta(E) q^{-2n}v-Ev\in V$, 
in which the top degree of $K$ is $n-1$. By repeating this operation 
finitely many times, we will finally get an element $0\neq a\omega$
with $a\in \C^{\ast}$. This means that $V=U_{q}(f(K))v$ for any $0
\neq v\in V$. So $V$ is irreducible. Therefore, we are done with the
proof. \qed

\begin{thm}
Let $V$ be a $U_{q}(f(K))-$module which admits a central
 character. Assume that $w\in V$ is a Whittaker vector.
Then the submodule $U_{q}(f(K))w\subset V$ is irreducible.
\end{thm}
{\bf Proof:} First of all, $U_{q}(f(K))w$ is a Whittaker module. Since
$V$ has a central character, then $U_{q}(f(K))w$ has a central
character. Thus $U_{q}(f(K))w$ is an irreducible Whittaker module.
\qed

\begin{thm}
Let $V_{1},\,V_{2}$ be any two irreducible $U_{q}(f(K))-$modules 
with the same central character. If $V_{1}$ and $V_{2}$ 
contain Whittaker vectors, then these vectors are unique 
up to scalars. And furthermore $V_{1}$ and $V_{2}$ are 
isomorphic to each other as $U_{q}(f(K))-$modules.
\end{thm}
{\bf Proof:} Since $V_{i}$ are irreducible and have Whittaker vectors,
then they are irreducible Whittaker modules. In addition, they have a 
central character, so the subspace of Whittaker vectors is one dimensional,
hence the Whittaker vectors are unique up to scalars. In this case, it
is obvious that they are isomorphic to each other. \qed

\subsection{The submodule structure of  a Whittaker module $(V,\,w)$}
In this section, we spell out the details 
about the structure of submodules of a 
Whittaker module $(V,\,w)$. We have a clean 
description about all submodules through 
the algebraic geometry of the affine line
$\mathbb{A}^{1}$. Throughout this section, 
we fix a Whittaker module $(V,w)$ of type 
$\eta$ and a cyclic vector $w$ of $V$. 
Our argument is more or less the same 
as the one in \cite{O}.
\begin{lem}
Let $Z(U_{q}(f(K)))=\C[\Omega]$ be the center 
of $U_{q}(f(K))$, then any maximal ideal 
of $Z(U_{q}(f(K)))$ is of the form $(\Omega-a)$ 
for some $a\in \C$. 
\end{lem}
{\bf Proof:} This fact follows directly from the fact that $\C$ 
is an algebraically closed field and Hilbert Nullstenllenzatz Theorem.\qed

Let $Z_{V}$ be the annihilator of $V$ in $Z(U_{q}(f(K)))$, then $Z_{V}=(g(\Omega))$ is the two sided ideal of $\C[\Omega]$ generated by some polynomial $g(\Omega)\in Z(U_{q}(f(K)))$. Suppose that we have a decomposition for $g$ as follows:
\[
g=\prod_{i=1,\,2\cdots,\,m} g_{i}^{n_{i}},
\]
where $g_{i}$ are irreducible. Then we have the following:
\begin{prop}
$V_{i}=U_{q}(f(K))\prod_{j\neq i} g_{j}^{n_{j}}w$ 
are indecomposable submodules of $V$. 
In particular we have
\[
V=V_{1}\oplus \cdots \oplus V_{m}
\] 
as a direct sum of submodules.
\end{prop}
{\bf Proof:} It is easy to verify that $V_{i}$ 
are submodules. Now we show each $V_{i}$ is 
indecomposable. Suppose not, we can assume 
without loss of generality that $V_{1}=W_{1}\oplus W_{2}$. 
Now $ Z_{V}=Z_{W_{1}}\cap Z_{W_{2}}$. Since $Z(U_{q}(f(K)))$ 
is a Principal Ideal Domain, hence $Z_{W_{i}}=(g_{i}(\Omega))$.
Thus we have $g_{i}\mid f_{1}^{n_{1}}$. This implies 
that the decomposition is not a direct sum. Therefore
$V_{i}$ are all indecomposable. Finally, the decomposition 
follows from the Chinese Reminder Theorem.\qed
\begin{prop}
Let $(V,\,w)$ be a Whittaker module and 
$Z_{V}=<g^{n}>$ where $g$ is an irreducible 
polynomial in $\C[\Omega]$. Let $V_{i}=U_{q}(f(K))g^{i}w, 
i=0,\,\cdots,\,n$ and $S_{i}=V_{i}/V_{i+1}, 
i=0,\,\cdots,\,n-1$. Then $S_{i},\,i=0,\,\cdots,\, n-1$ 
are irreducible Whittaker modules of 
the same type $\eta$ and form a composition 
series of $V$. In particular $V$ is of 
finite length. 
\end{prop}
{\bf Proof:} The proof follows from the fact that $Z_{S_{i}}=<g>$ for all $i$.\qed

With the same assumption, we have the following
\begin{cor}
$V$ has a unique maximal submodule $V_{1}$.
\end{cor}
{\bf Proof:} This is obvious because the only 
maximal ideal of $Z_{V}$ is $<g>$.
\qed

Based on the above propositions, the irreducibility and indecomposability 
are reduced to the structure of $Z_{V}$. $V$ is irreducible if and only if $Z_{V}$ 
is maximal in $Z(U_{q}(f(K)))$. And $V$ is indecomposable if and only if $Z_{V}$ 
is a primary ideal of $Z(U_{q}(f(K)))$. The following proposition is a refinement 
of the submodule structure of $(V,\,w)$.
\begin{prop}
Suppose $(V,\,w)$ is an indecomposable Whittaker module with $Z_{V}=<g^{n}>$, 
then any submodule $T\subset V$ is of the form
\[
T=U_{q}(f(K))g^{i}w
\] 
for some $i\in \{0,\,\cdots,\,n\}$.
\end{prop}
{\bf Proof:} The proof is obvious.
\qed

Now we investigate the submodule structure of any Whittaker module 
$(V,\,w)$ with a nontrivial central annihilator $Z_{V}$. First of all, 
we recall some notations from \cite{K}. Let $T\subset V$ be any 
submodule of $V$, we define an ideal of $Z$ as follows:
\[
Z(T)=\{x\in Z\mid xT\subset T\}
\]
We may call $Z(T)$ the normalizer of $T$ in $Z$. Conversely for any ideal $J\subset Z$ 
containing $Z_{V}$, $JV\subset V$ is a submodule of $V$. We have the following Theorem:
\begin{thm}
Let $(V,w)$ be a Whittaker module with $Z_{V}\neq 0$. Then there is 
a one-to-one correspondence between the set of all submodules 
of $V$ and the set of all ideals of $Z$ containing $Z_{V}$ given by the 
\[
T\longrightarrow Z(T)\]
 and 
\[
J\longrightarrow JV
\]
which are inverse to each other.
\end{thm} 
{\bf Proof:} The proof is a straightforward.
\qed

Now we give a description of the basis of any Whittaker module $(V,w)$.
\begin{prop}
Let $(V,w)$ be a Whittaker module and suppose that $Z_{V}=<g(\Omega)>$ where $g\neq 0$ 
is monic polynomial of degree $n$. Then 
\[
{\bf B}=\{F^{i}K^{j}\mid 0\neq i\leq n-1, j\in {\Z}\}
\]
is a $\C-$basis of $V$. If $g=0$, then 
\[
{\bf B}=\{F^{i}K^{j}\mid i\geq 0, j\in {\Z}\}
\]
is a $\C-$basis of $V$.
\end{prop}\qed

\subsection{The annihilator of a Whittaker module}
In \cite{O}, it was proved the annihilator of any Whittaker module
of $U_{q}(sl_{2})$ is centrally generated. In this section we
generalize this result to our situation. We closely follow the 
approach in \cite{O}. First of all, we need some lemmas:
\begin{lem}
Let $(V,\,w)$ be a Whittaker module of type $\eta$ with a fixed
Whittaker vector $w$. Suppose there is a $u\in U_{q}(f(K))$ such that
$uK^{i}w=0$ for all $i>0$. Then $uw=0$.
\end{lem}
{\bf Proof:} (We will adopt the proof of Lemma 6.1 in \cite{O}). We 
can write $u=\sum_{i\in \Z}x_{i}$ where $Kx_{i}=q^{i}x_{i}K$. Suppose 
the statement is false, then there exists a minimal such $u$ with 
respect to the length of the above expression of $u$. We may assume
$u$ has more than one summand, otherwise $0=x_{r}Kw=q^{r}Kx_{r}w$
implies that $uw=0$. Since we have $uK^{i}w=$ for all $i>0$. In
particular, we have 
\[
0=(\sum_{r}x_{r})Kw
\]  
and 
\[
0=(\sum_{r}x_{r})K^{2}w=K(\sum_{r}q^{-r}x_{r})Kw.
\]
Thus we have $0=(\sum_{r}q^{-r}x_{r})Kw$. Hence both $u=\sum_{r}x_{r}$
and $u'=\sum_{r}q^{-r}x_{r}$ annihilate $Kw$. Note that
\[
Ku'K^{i}w=K(\sum_{r}x_{r})K^{i}w=(\sum_{r}x_{r})K^{i+1}w
=uK^{i+1}w=0\]
for any $i>0$. Thus $0=K^{-1}Ku'K^{i}w=u'K^{i}w$. Now let
$m=max\{r\mid x_{r}\neq 0\}$, and note that $u-q^{m}u'$ annihilates
$K^{i}w$ for $i>0$. But 
\[
u-q^{m}u'=u-q^{m}(\sum_{r}q^{-r}x_{r})=\sum_{r\neq m}(1-q^{m-r})x_{r}
\]
has fewer nonzero terms than $u$. This is a contradiction.\qed

The following corollary can be proved the same as in \cite{O}:
\begin{cor}
Let $V$ be a Whittaker module with $Z_{V}=0$, then $U_{V}=UZ_{V}=0$.
\end{cor}\qed

In addition, the following similar lemma holds:
\begin{lem}
Assume that for any simple Whittaker module $V$, we have
$U_{V}=UZ_{V}$. Then for any Whittaker module $V$ of finite length, $U_{V}=UZ_{V}$.
\end{lem}\qed

As remarked in \cite{O}, the above lemma reduces the problem to the
case where $V$ is an irreducible Whittaker module with a nonzero
central annihilator $Z_{V}=<\Omega-a>$ for some $a\in \C$. To
summarize the result, we have the following:
\begin{thm}
Let $V$ be a Whittaker module, then $Ann(V)$ is centrally generated, i.e. 
\[
Ann_{U_{q}(f(K))}(V)=U_{q}(f(K))Ann_{Z(U_{q}(f(K)))}(V)
\]
\end{thm}
{\bf Proof:} (This argument is essentially due to Smith \cite{S}). First of all, the primitive ideal $U_{V}$ has infinite codimension
since $U_{V}\subset U_{w}$ and $U/U_{w}=V$ is infinite-dimensional. 
Since $U_{q}(f(K))$ is a domain, then $(0)$ is a
prime ideal of $U_{q}(f(K))$. $U_{q}(f(K))Z_{V}$ is also prime. Since
$U_{q}(f(K))_{V}$ is primitive, it is prime. Thus we have a chain of
prime ideals: $(0)\subset U_{q}(f(K))Z_{V}\subset U_{q}(f(K))_{V}$. 
Let $R=U/U_{V}$, then $R$ is a primitive ring. If $R$ is artinian,
then $R$ is finite dimensional which is contradicting to the fact that
$U_{V}$ has infinite codimension. So $R$ is not artinian. Now $U$ has
GK-dimension $3$. Suppose that the prime ideals $(0)\subset
UZ_{V}\subset U_{V}$ are different, then GK-dimension of $U/U_{V}$ is
at most $1$. But there are no non-artinian finitely generated 
noetherian primitive $\C-$algebras with GK-dimension $1$. Thus $R$ has
to be finite dimensional, which is a contradiction. So we are done with
proof.\qed

\end{document}